\documentclass[12pt]{amsart}
\usepackage{amsmath}
\usepackage{amsfonts,amssymb,amsthm, txfonts, pxfonts,amscd}
\def\struckint{\mathop{%
\def\mathpalette##1##2{\mathchoice{##1\displaystyle##2}%
  {##1\textstyle##2}{##1\scriptstyle##2}{##1\scriptscriptstyle##2}}%
\mathpalette
{\vbox\bgroup\baselineskip0pt\lineskiplimit-1000pt\lineskip-1000pt
\halign\bgroup\hfill$}
{##$\hfill\cr{\intop}\cr\diagup\cr\egroup\egroup}%
}\limits}

\newcommand{\specrad}[1]{\mathcal{R}(#1)}

\newtheorem{theorem}{Theorem}[section]
\newtheorem{lemma}[theorem]{Lemma}
\newtheorem{corollary}[theorem]{Corollary}

\newtheorem{definition}[theorem]{Definition}
\newtheorem{fact}[theorem]{Fact}

\newtheorem{notation}[theorem]{Notation}
\newtheorem{theorem-definition}[theorem]{Theorem-Definition}

\theoremstyle{remark}

\newtheorem{remark}[theorem]{Remark}
\newtheorem{condition}[theorem]{Condition}

\newtheorem{example}[theorem]{Example}

\newcommand{\cx}{\mathbb{C}}
\newcommand{\integers}{\mathbb{Z}}

\newcommand{\reals}{\mathbb{R}}

\DeclareMathOperator{\diag}{Diag}

\DeclareMathOperator{\tr}{tr}

\DeclareMathOperator{\Sp}{Sp}
\DeclareMathOperator{\SL}{SL}
\DeclareMathOperator{\GL}{GL}
\DeclareMathOperator{\Gal}{Gal}

\providecommand{\matnorm}[1]{\lvert\lVert#1\rVert\rvert}
\providecommand{\opnorm}[1]{\lvert\lVert#1\rVert\rvert_{\text{op}}}
\providecommand{\norm}[1]{\lVert#1\rVert}
\providecommand{\ffp}[1]{\mathbb{F}_{#1}}
\providecommand{\abs}[1]{\lvert#1\rvert}

\begin{document}

%-------------- Author entries --------------------

\title[Groups, polynomials, matrices, and
automorphisms]{Walks on groups, counting reducible matrices, polynomials, and surface and free group automorphisms}

\author{Igor Rivin}

\address{Department of Mathematics, Temple University, Philadelphia}

\email{rivin@math.temple.edu}

\thanks{I  would like to thank Ilya Kapovich for asking the
  questions on the irreducibility of  automorphisms of surfaces and
  free groups, and for suggesting that these questions may be
  fruitfully attacked by studying the action on homology.
  I would also like to thank Nick Katz for making him aware of Nick
Chavdarov's work, Sinai Robins for making him aware of Morris Newman's
classic book, and Akshay Venkatesh, Peter Sarnak, and Alex Lubotzky for enlightening
conversations. The author would also like to thank Benson Farb, Ilya Kapovich, and Lee Mosher for
comments on previous versions of this note.}

\curraddr{Mathematics Department, Stanford University, Stanford, California}
\email{rivin@math.stanford.edu}

\date{\today}


%\subjclass{57M50, 32G15}

\begin{abstract}
We prove sharp limit theorems on random walks on graphs with values in
finite groups. We then apply these results (together with some
elementary algebraic geometry, number theory, and representation
theory) to finite 
quotients of lattices in semisimple Lie groups (specifically
$\SL(n,\integers)$ and $\Sp(2n, \integers)$) to show that a ``random''
element in one of these lattices has irreducible characteristic
polynomials (over $\integers$). The term ``random'' can be defined in
at least two ways (in terms of height and also in terms of word length
in terms of a generating set) -- we show the result using both
definitions. 

We use the above results to show that a random (in terms of word
length) element of the mapping class group of a surface is
pseudo-Anosov, and that a a random free group automorphism is
irreducible with irreducible powers (or \emph{strongly irreducible}
\footnote{terminology due to Lee Mosher}).
\end{abstract}

\maketitle

\section*{Introduction} 

This paper was inspired by the following question, first brought to
the author's attention by Ilya Kapovich:

\begin{quotation}
Is it true that a random element of the mapping class group of a
surface is pseudo-Anosov?
\end{quotation}

The definition of \emph{random} in the question above is not
explicitely given, but a reasonable way to define is to fix a
generating set of the mapping class group, and look at all the words
of bounded length.\footnote{Another way is to look at a combinatorial
  ball of radius $N$ around identity; it turns out that our results,
  together with the results announced recently by Ursula Hammenstadt
  answer both questions.}
Kapovich had suggested that a reasonable way to attach this question
was to study the action of the mapping class group $\mathcal{M}_g$ on
homology, which gives a symplectic representation of $\mathcal{M}_g.$
Results of Casson and Bleiler (see \cite{CassonBleiler}) then give a
set of sufficient conditions for a surface automorphism to be
pseudo-Anosov in terms of its image under that representation.

This paper, then, is the embodiment of the program as described
above. Along the way, we show a number of results not explicitely
related to low dimensional topology or geometric group theory.  

Here is a summary of the results:

In my old preprint \cite{walks} I state, and sketch the proof of a
general equidistribution theorem for products of elements of finite
(and, more generally, compact) groups along long paths in finite
graphs. In the current paper I give complete arguments. In the follow
up-paper \cite{rirred2} much sharper results are given
(with speed of convergence bounds). In particular, the
convergence bounds are \emph{uniform} for families of finite quotients
of a group satisfying Lubotzky's property $\tau$ (or the stronger
Kazhdan's property $T$). 

We then show that the set of polynomials reducible
over $\integers$ with constant coefficient $1$ is an algebraic
subvariety of the set of polynomials, and using elementary theory of
algebraic groups show that the set of matrices in $\SL(n, R)$ (where
$R$ is the coefficient ring) with reducible characteristic polynomials
is a finite union of Zarisky-closed sets, so in particular, when
$R=\integers_p,$ the proportion of matrices in $\SL(n, R)$ with
reducible characteristic polynomials decreases as $O(1/p),$ with
constants independent of $p$ (Section \ref{matrices}). Using these
methods together with Weil's estimate for the number of points on a
curve over $\mathbb{F}_p$ we show that the probability that the
characteristic polynomial of a matrix in $\SL(n, \integers_p)$ has the
full symmetric group $S_n$ as Galois group goes is bounded (uniformly)
away from $0.$ Along the way, we show that the probability that a
polynomial with fixed constant term $1$ has a prescribed splitting
type $T$ modulo $p$ is the same as the probability that a
\emph{random} polynomial has splitting type $T$ (for large
$p$). (Section \ref{galoisrestricted}). 

We quote results of Borel (proofs can be found in Nick Chavdarov's
thesis \cite{chavdarov}) which show that the probability that a matrix
in 
$\Sp(2n, \integers_p)$ has reducible characteristic polynomial is
bounded away from $1$ (uniformly in $p.$) (Section \ref{matrices2})

Combining these results with with the results on walks on graphs,
together with elementary local-to-global estimates, we show that for a given
\emph{undirected} graph $G$ whose vertices are decorated with a
symmetric generating set of $\Gamma$ (where $\Gamma$ is either $\SL(n,
\integers)$ or $\Sp(2n, \integers)$), the probability that the product
of generators along a random long walk of length $N$ gives a matrix
with reducible characteristic polynomial (or, in the case of $\SL(n,
\integers),$ a polynomial with non-generic Galois group) goes to $0$
as $N$ goes to infinity. 

From a number-theoretic standpoint, using the length of the
representing word in generators is a less natural way to define the
size of an integral matrix than a(ny) matrix norm. It turns out that
with that definition of size the result follows by combining our
results on finite matrix groups with the work of \cite{DRS}, and
effective bounds can be obtained using
the uniformity results of \cite{NevoSarnak} (see Theorem \ref{slnzgro}).

Finally, we apply the above-mentioned results to automorphisms of
surfaces (Section \ref{mcg}) and free group automorphisms (Section
\ref{free}). For surface automorphisms we use the observations of
Casson and Bleiler (\cite{CassonBleiler}) to show that for an
arbitrary generating set, most products of up to $N$ generators are
pseudo-Anosov, and if the generating set happens to be symmetric, then
the fraction of the non-pseudo-Anosov products goes to zero
exponentially fast with $N.$ It can be argued that it is more natural
to consider all elements in a combinatorial ball of radius $N.$ Since
our limit theorems for graphs are for undirected graphs, such a result
will follow (for certain generating sets) once we know that the
mapping class groups are bi-automatic.

Similarly, it is shown that most words in a generating set of the
outer automorphism group of a free group $F_n$ are irreducible with
irreducible powers, or what Lee Mosher calls strongly irreducible.

\section{Generalities on algebraic geometry and algebraic groups}
\label{closure}
 
First, recall the following:
\begin{definition}[\cite{geck2004}[p.~23]
Let $V\subseteq k^n,$ and $W \subseteq k^m$ be non-empty algebraic
sets. We say that $\phi: V \rightarrow W$ is a \emph{regular} map, if
there exist $f_1, \dotsc, f_m \in k[X_1, \dotsc, X_n],$ (where $X_1,
\dotsc, X_n$ are indeterminantes) such that 
\[
\phi(x) = (f_1(x), \dotsc, f_m(x))
\]
for all $x \in V.$
\end{definition}

\begin{example}
\label{firsteg}
Let $V = \SL(n, k) \subset k^{n^2},$ and let $W = k^n.$ 
Then the map associating to each matrix the coefficients its characteristic polynomial
is a regular map.
\end{example}

Consider the following setup: we have a parametrized set $S$ in $k^n$
(for $k$ an algebraically closed field), that is:
\begin{gather*}
x_1 = f_1(s_1, \dotsc, s_m),\\
x_2 = f_2(s_1, \dotsc, s_m),\\
\vdots\\
x_n = f_n(s_1, \dotsc, s_m),
\end{gather*}
where $f_1, \dotsc, f_n$ are polynomials in $s_1, \dotsc, s_m.$
By the ``Implicitization Theorem'' \cite{CLO_IVA}[Chapter 3], the
Zariski closure of $S$ is an affine variety, whose dimension is
bounded above by $m$ (by considering the dimension of the tangent
space, see \cite{geck2004}[Theorem 1.4.11]).
\subsection{Fibers of dominant morphisms}
\label{domfib}
The following results are standard; the statement is taken from
\cite{geck2004}[pp.~116,118]. First, let $A[X]$ be the algebra of
regular functions on $X,$ and second, for $g \in A[X],$ let $X_f =
\{x\in X \vert f(x) \neq 0\}.$ Now
\begin{theorem}
\label{geck1}
Let $\phi: X \rightarrow Y$ be a dominant morphism of irreducible
affine varieties. In particular, $\phi^*:A[Y] \rightarrow A[X]$ is
injective, and $d = \dim X - \dim Y \geq 0.$ Then we have a
factorization
\[
\begin{CD}
\phi\vert_{X_{\phi^*(g)}} : X_{\phi^*(g)}
  @>{\overline{\phi}}>> Y_g \times k^d @>{p_1}>> Y_g,
\end{CD}
\]
where $\overline{\phi}$ is a finite dominant morphism, and $p_1$ is
the projection on the first factor.
\end{theorem}

\begin{corollary}
\label{borel}
Let $\alpha: X \rightarrow Y$ be a dominant morphism of irreducible
varieties, and put $r = \dim X - \dim Y.$ Then, there exists a
non-empty open set $U \subseteq Y,$ such that $U \subseteq \phi(X),$
and such that $\dim \phi^{-1}(y) = r$ for all $y \in U.$ 
\end{corollary}

\begin{remark}
\label{geckrem}
$U$ can be taken to be of the form $X_f,$ for some $f \in A[Y].$ 
\end{remark}
A \emph{dominant morphism $\alpha$} is a morphism such that the
preimage of any open dense set in the image variety is dense. In
particular, a surjective regular map is dominant.

\begin{example}
\label{secondeg}
Consider $V\subset k^{n^2},$ where $V = \SL(n, k).$  The map $\chi,$
which associates to each matrix the coefficients of its characteristic
polynomial is a dominant morphism onto $k^{n-1}$ (by the companion matrix
construction, it is surjective). By the results (and notation) above,
there is a $g \in A[V],$ such that the fibers of $\chi$ restricted to
$V_g$ are $n^2-n$-dimensional varieties.

Similarly, let $\Sp(2n, k) = W \subset k^{4n^2},$  and let
$\chi^\prime$ be the map which associates to each matrix $m$ the
coefficients of $x, \dotsc, x^n$ of the characteristic polynomial of
$m.$
By Theorem \ref{kirby}, $\chi^\prime$ is a dominant morphism, and by
Fact \ref{dickfact} and the results of this section, there is an $h
\in A[W],$ such that the fibers of $\chi^\prime$ restricted to $W_h$
are $2n^2$-dimensional varieties.
\end{example}

\section{Counting points on varieties}
\label{counting}

Let $S$ be a variety o
f dimension $m$ over $\mathbb{C}$ (or $\overline{\mathbb{Q}}.$  Consider a reduction of $S$
modulo $p.$

\begin{theorem}[Lang-Weil, \cite{langweil}]
\label{langweil}
The number of $F_p$ points on $S$ grows as $O(p^m).$ The implied
constant is uniform (that is, it is a function of the dimension and
codimension of the variety \emph{only}).
\end{theorem}

It should be noted that we are using this theorem for the upper bound
only, and the upper bound is an easy result (see \cite[Lemma
4.1.3]{geck2004}), unlike the full Lang-Weil result which gives both
upper and lower bounds. Lang and Weil is deduced from:
A.~Weil's estimate on the number of
$\mathbb{F}_p$ points on a curve defined over $\mathbb{F}_p:$
\begin{theorem}[A. Weil,\cite{weilcourbes}]
\label{weilcourbesthm}
Let $f \in \mathbb{F}_p[X, y]$ be an absolutely irreducible (that is,
irreducible in $\overline{\mathbb{F}}_p[X, Y]$) polynomial of degree $d$. Then if
\[\mathcal{C} = \{(x, y) \in \mathbb{F}_p^2 | f(x, y) = 0\},\]
we have the estimate
\[||C| - p| \leq 2 g \sqrt{p} + d^2,\] where $g$ is the genus of the
curve defined by $f$ (which satisfies $g \leq (d-1)(d-2).$ This estimate is sharp.
\end{theorem}

Theorem \ref{langweil}, together with the results of Section
\ref{domfib} imply:
\begin{theorem}
\label{domcount}
With the situation as described in Corollary \ref{borel}, if the
ground field is $\ffp{p},$ if $S \subseteq Y,$ then $|\alpha^{-1}(S)| \leq
c_1 |S| p^r+c_2(p^{\dim X - 1})$ where the constants $c_1, c_2$
depend only on the dimensions of $X, Y.$
\end{theorem}

\section{Classical groups}
\label{classgroups}
In this paper we will be primarily concerned with the special linear
group and the symplectic group over various domains. We will need the
following facts:

\begin{fact}[See \cite{borelgroups,geck2004}]
\label{irredgp}
The groups $\GL(n, k), \SL(n, k)$ and $\Sp(n, k)$ over an algebraically closed
field $k$ are \emph{irreducible} as algebraic varieties. This is
equivalent to saying that the groups in question are connected.
\end{fact}
\begin{fact}
\label{dimgp}
The dimension of $\SL(n, k)$ equals $n^2-1.$ The dimension of $\Sp(2n,
k)$ equals $2n^2 + n.$ In both cases, dimension is meant in the sense
of algebraic geometry.
\end{fact}

Let $\ffp{p} = \mathbb{Z}/p \mathbb{Z}.$

The following goes back to Dickson:
\begin{fact}
\label{dickfact}
The order of $\SL(n, \ffp{p})$ equals 
\[
p^{(n^2-n)/2} (p^2-1)(p^3-1) \dots (p^n - 1) = p^{n^2-1} +
O(p^{n^2-3}).
\]
The order of $\GL(n, \ffp{p})$ equals
\[
p^{(n^2-n)/2} (p-1)(p^2-1)(p^3-1) \dots (p^n - 1) = p^{n^2} +
O(p^{n^2-1}).
\]
The order of $\Sp(2n, \ffp{p})$ equals
\[
p^{n^2}(p^2-1)(p^4 - 1) \dots (p^{2n} - 1) = p^{2n^2+n} + O(p^{2n^2 +
  n - 2}).
\]
\end{fact}
 \section{Applications to polynomials}
\label{poly}
 Let $\mathcal{P}_d(k, \beta)$ be the set of all monic  polynomials in one
variable of 
 degree $d$ over a field 
 $k$ with constant coefficient $\beta,$ and let $\mathcal{P}_d^m(k,
 \beta, \alpha)\subseteq \mathcal{P}_d(k, \beta)$ be the
 set of those polynomials which have a polynomial factor (over
 $k$) of degree $m$ with constant term $\alpha.$  
Let us identify $\mathcal{P}_d(k, \beta)$ with the affine space $k^{d-1}.$
Then, we have the following: 
\begin{theorem}
    \label{hyper1}
    The set $\mathcal{P}_d^m(k, \beta, \alpha)$ is contained in an affine hypersurface of 
    $\mathcal{P}_{d}(k, \beta).$
\end{theorem}

\begin{proof}
    Let \[p(x) = x^{d} + \sum_{i=0}^{d-1}a_{i} x^{i}¥ \in 
    \mathcal{P}.\] By assumption, $p(x) = q(x) r(x).$ Assume that the 
    degree of $q(x) = m,$ while the constant term of $q(x)$ equals 
    $\alpha.$ Writing \[q(x) = x^{m} + \sum_{j=1}^{m-1} b_{j}x^{j},\] 
    and \[r(x) = x^{d-m} + a_{0}/\alpha + \sum_{k=1}^{d-m-1}
c_{k}x^{k},\] we see that $\mathcal{P}_d^m(k, \beta,\alpha)$ is a polynomially
parametrized hypersurface of $\mathcal{P}_d(\beta,k).$
\end{proof}

\begin{corollary}
\label{allreds}
The set of all reducible monic polynomials with fixed constant term
$\beta$ is the union 
\[
\bigcup_{m=1}^d \mathcal{P}_d^m(k, \beta, \alpha)
\]
of affine hypersurfaces of $\mathcal{P}_d(k, \beta) \simeq k^{d-1}.$
\end{corollary}

\begin{theorem}
\label{constone2}
Let $P(d)(\ffp{p})$ be the set of polynomials of degree $d$ with
coefficieded in absolute value by $N$ and constant coefficient $1,$ and let
$R_1(d)(\ffp{p})$ be the set of polynomials reducible over $\ffp{p}$ with
some factor having constant term equal to $1 \mod p.$ Then, $R_1(d)$ lies
on an algebraic 
hypersurface $\mathbb{C}^{d-1}$ (where the coordinates are the
coefficients), and consequently
\[
\dfrac{R_1(d)}{P_1(d)} = O\left(\dfrac{1}{p}\right).
\]
\end{theorem}

\begin{proof}
Immediate from Theorem \ref{langweil}.
\end{proof}

Finally, one definition:
\begin{definition}
We say that a polynomial $p(x) \in P_1(d)$ is  \emph{reciprocal} if 
$x^d p(1/x) = p(x)$ -- in other words, the list of coefficients of $p$
is the same read from left to right as from right to left. Reciprocal
polynomials can also be defined as follows:
A (monic) polynomial (of even degree $2n$) is reciprocal if it can be
written as
\[
\prod_{j=1}^n (x-r_i)(x-1/r_i) = \prod_{j=1}^n (x^2 - (r_i + r_i^{-1})
x + 1).
\]
\end{definition}

\section{Applications to matrices}
\label{matrices}

\subsection{The special linear group.}
\label{slnsec}
\begin{lemma}
\label{slnpprob}
The probability that the characteristic polynomial of a matrix $M$ in
$\SL(n, \ffp{p})$ has a factor with constant term $1$ is of order $O(1/p).$
\end{lemma}
\begin{proof}
This follows immediately from Example \ref{secondeg} and Theorems
\ref{constone2} and \ref{domcount}
\end{proof}

Unfortunately, since the number of integral points on $\SL(n,
\mathbb{Z})$ of height (absolute value) bounded by $B$ grows much slower than
$B^{n^2-1}$ the above results \emph{do not} imply the following
\begin{theorem}
\label{slnzgro}
The probability that a matrix in $\SL(n, \mathbb{Z})$ with coefficients
bounded by $B$ has reducible characteristic polynomial goes to $0$ as
$B$ goes to infinity.
\end{theorem}
\begin{proof}
For the ``soft'' result as above we use results of Duke Rudnick and
Sarnak (\cite{DRS}) or of Eskin and McMullen (\cite{EskinMac}). These
imply that the matrices in $\SL(n, \mathbb{Z})$ with coefficients
bounded by $B$ are asymptotically equidistributed among the cosets of
any subgroup of finite index, in particular the cosets of the principal
congruence subgroup $\Gamma_p(n, \mathbb{Z}).$ Since reducibility
modulo $p$ depends only on the coset modulo $\Gamma_p,$ Lemma
\ref{slnpprob} immediately replies the result.
\end{proof}
\begin{remark}
To get an effective result, we need to know that the error terms in
equidistribution modulo $\Gamma_p$ are uniform (do not depend on
$p.$). Precisely such a result is shown in the preprint of Amos Nevo
and Peter Sarnak \cite{NevoSarnak}.
\end{remark}

\section{Random products of matrices in the symplectic and special
linear groups}
\label{matrices2}

In the preceeding section we defined the size of a matrix by (in
essence) its $L^1$ norm (any other Banach norm will give the same
results). However, it is sometimes more natural to measure size
differently: In particular, if we have a generating set $\gamma_1,
\dots, \gamma_l$ of our lattice $\Gamma$ (which might be $\SL(n,
\mathbb{Z})$ or $Sp(2 n, \mathbb{Z})$) we might want to measure the
size of an element by the length of the (shortest) word in $\gamma_i$
equal to that element -- this is the combinatorial measure of
size. The relationship between the size of elements and combinatorial
length is not at all clear, so the results in this section are proved
quite differently from the results in the preceding section. 

We will be using Theorems \ref{mythm1} and \ref{mythm2}.

\begin{remark}
We will be applying Theorem \ref{mythm1} to groups
$\SL(n, \mathbb{Z}/p \mathbb{Z})$ and $\Sp(2n, \mathbb{Z}/p
\mathbb{Z}).$ Since those groups have no non-trivial one-dimensional
representations, the assumption on $\rho$ in the statement of the
theorem is vacuous.
\end{remark}

We will also need the following results of Nick Chavdarov and Armand Borel.

\begin{theorem}[N.~Chavdarov, A.~Borel \cite{chavdarov}]
\label{chavthm}
Let $q > 4$, and let $R_q(n)$ be the set of $2n \times 2n$ symplectic matrices over the field $F_q$ with
\emph{reducible} characteristic polynomials. Then 
\[
\dfrac{|R_q(n)|}{|\Sp(2n, F_p)|} < 1- \frac{1}{3 n}.
\]
\end{theorem}

\begin{theorem}[N.~Chavdarov, A.~Borel \cite{chavdarov}]
\label{chavthm2}
Let $q > 4$, and let $G_q(n)$ be the set of $n \times n$ matrices with
determinant $\gamma \neq 0$ over the field $F_q$ with
\emph{reducible} characteristic polynomials. Then 
\[
\dfrac{|G_q(n)|}{|\SL(n, F_q)|} < 1- \frac{1}{2 n}.
\]
\end{theorem}

Theorem \ref{chavthm2} follows easily from the following result of A.~Borel:
\begin{theorem}[A.~Borel]
\label{borthm2}
Let $F$ be a monic polynomial of degree $N$ over $\mathbb{Z}/p \mathbb{Z}$
with nonzero constant term.  Then, the number $\#{F, p}$ of matrices in $\GL(N,
p)$ with characteristic polynomial equal to $F$ satisfies
\[
(p-3)^{N^2-N} \leq \#(F, p) \leq (p+3)^{N^2-N}.
\]
\end{theorem}
Theorem \ref{borthm2} will be used in Section \ref{strong}.

We now have our results:

\begin{theorem} 
\label{randprodthm}
Let $G$ and $S_N,$ be as in the statement of Theorem
\ref{mythm1}, but with $\Gamma = \Sp(2n, \mathbb{Z}),$ or $\Gamma = \SL(2,
\mathbb{Z}).$ 
 Then the
probability that a matrix in $S_N$ has a reducible characteristic
polynomial goes to $0$ as $N$ tends to infinity.
\end{theorem}
\begin{proof}  Let $\Gamma_l$ be the set of matrices in $\Gamma$ reduced
  modulo $l$ -- it is known (see \cite{newmanmats}) that $\Gamma_l$ is $\SL(n,
  l)$ or $\Sp(2n, l)$ (depending on which $\Gamma$ we took. Let $p_1, \dotsc,
  p_k$ be distinct primes, let $K = p_1 
\dots p_k.$ 
We know that:
\[\Gamma_K = \Gamma_{p_1} \times \dots
\times \Gamma_{p_k}.\] (see \cite{newmanmats} for the proof of the last
equality). A generating set of $\Sp(2n, \mathbb{Z})$ projects via
reduction modulo $K$ to a a generating set of $\Gamma_K$ (see, again,
Newman's book \cite{newmanmats}), and also, via reduction mod $p_i$ to
generating sets of the $\Sp(2n, p_i).$ By Theorems
\ref{mythm1},\ref{mythm2} and 
\ref{chavthm}, the probability that the characteristic polynomial in a
random product of $N \gg 1$ generators is reducible modulo all of the
$p_i$ is at most equal to $(1-3/n)^k.$ Since this is an upper bound on
the probability of being reducible modulo $\mathbb{Z},$ the result
follows.
\end{proof}

\begin{remark}
Using Lemma \ref{slnpprob} instead of Theorem \ref{chavthm2} for
$\SL(n, \mathbb{Z})$ gives a sharper result, as well as a more elementary argument.
\end{remark}

\begin{remark}
The proof of Theorem \ref{slnzgro} translates \emph{mutatis mutandis}
into this setting -- instead of the principal congruence subgroup
$\Gamma_p$ for a single prime $p,$ we use it for $q = p_1 \times \dots
\times p_k.$
\end{remark}

An example of a graph $G$ is a graph where every vertex is connected
to every vertex (\emph{including itself}), and the set of label is a
\emph{symmetric} generating set.
In this case, we are just taking
random products of generators or their inverses. Another is the graph
(studied  in \cite{walks}) where a vertex labelled by a generator $a$
is connected to every vertex \emph{except} the one labelled by $a^{-1},$
so that only reduced words in the generators are allowed, and so on.

\section{Stronger irreducibility}
\label{strong}
We might ask if something stronger than irreducibility of the
characteristic polynomial can be shown.The answer is in the
affirmative. Indeed, the methods of the preceeding sections combined with the
results of the Appendix give immediately:
\begin{theorem}
\label{galthm}
The probability that a random word of length $L$ in a generating set of $SL(N,
\mathbb{Z})$ has characteristic polynomial with Galois group $S_N$
goes to $1$ as $L$ goes to infinity.
\end{theorem}

Aside from its intrinsic interest, Theorem \ref{galthm} implies the
following:
\begin{theorem}
\label{iwip}
The probability that a random word $w$ of length $L$ in a generating set
of $\SL(n, \mathbb{Z})$ and all proper powers $w^k$ have irreducible
characteristic polynomials goes to $1$ as $L$ goes to infinity.
\end{theorem}

Theorem \ref{iwip} will follow easily from Theorem \ref{galthm}
together with the following Lemma\footnote{Compare with
  \cite{chavdarov}[Lemma 5.3]}:
\begin{lemma}
\label{iwip2}
Let $M \in \SL(n, \mathbb{Z})$ be such that the characteristic
polynomial of $M^k$ is \emph{reducible} for some $k.$ Then the Galois group
of the characteristic polynomial of $M$ is \emph{imprimitive}, or the characteristic polynomial of $M$
is cyclotomic.
\end{lemma}
\begin{remark} For the definition of \emph{imprimitive} see, for
  example, \cite{wielandt,mhall}.
\end{remark}
\begin{proof}
Assume that the characteristic polynomial $\chi(M)$ is irreducible
(otherwise the conclusion of the Lemma obviously holds, since the
Galois group of $\chi(M)$ is not even transitive). Let the roots of
$\chi(M)$ (in the algebraic closure of $\mathbb{Q}$) be $\alpha_1, \dotsc, \alpha_n.$ The roots of
$\chi(M^k)$ are $\beta_1, \dotsc, \beta_n,$ where $\beta_j = \alpha_j^k.$ Suppose that
$\chi(M^k)$ is reducible, and so there is a factor of $\chi(M^k)$
whose roots are $\beta_1, \dots, \beta_l,$ for some $l < n.$ Since
$\Gal(\chi(M))$ acts transitively on $\alpha_1, \dotsc, \alpha_n,$ it
must be true that for every $i \in \{1, \dots, n\},$ $\alpha_i^k =
\beta_j,$ for some $j \in \{1, \dotsc, l\}.$ Let $B_j$ be those $i$
for which $\alpha_i^k = \beta_j.$ This defines a partition of $\{1,
\dotsc, n\}$ into blocks, which is stabilized by the Galois group of
$\chi(M),$ and so $G$ is an intransitive subgroup of $S_n,$
\emph{unless} $l = 1.$ In that case, the characteristic polynomial of
$M^k$ equals $(x-\beta)^n,$ and since $M^k \in \SL(n, \mathbb{Z})$ it
follows that $\beta = 1,$ and all the eigenvalues of $M$ are $n$-th
roots of unity, so that $M^k = 1.$
\end{proof}

\section{The mapping class group} 
\label{mcg}
Let $S_g$ be a closed surface of genus $g\geq 1,$
 and let $\Gamma_g$ be the mapping class group of $S_g.$ The group $\Gamma_g$
 admits a homomorphism $\mathfrak{s}$ onto $\Sp(2g, \mathbb{Z})$  (we associate to each
 element its action on homology; the symplectic structure comes from the
 intersection pairing). The following result can be found in \cite{CassonBleiler}:
\begin{theorem}
\label{cassonthm}
For $\gamma \in \Gamma_g$ to be pseudo-Anosov, it is sufficient that $g =
\mathfrak{\gamma}$ satisfy all of the following conditions:
\begin{enumerate}
\item
The characteristic  polynomial of $g$ is irreducible.
\item
The characteristic polynomial of $g$ is not cyclotomic.
\item
The characteristic polynomial of $g$ is not of the form $g = h(x^k),$ for some
$k>1.$ 
\end{enumerate}
\end{theorem}
The following is a corollary of our results on matrix groups:
\begin{theorem}
Let $g_1, \dots, g_k$ be a generating set of $\Sp(2n, \mathbb{Z}).$ The
probability that a random product of length $N$ of  $g_1, \dots, g_k$
satisfies the conditions of Theorem \ref{cassonthm} goes to $1$ as $N$ goes to
infinity. 
\end{theorem}
\begin{proof}
We prove that the probability that the random word $w_N$ \emph{not} satisfy the
conditions goes to $0.$ By Theorem \ref{randprodthm}, the probability that
$w_N$ has reducible characteristic polynomial goes to $0.$ In order for the
characteristic polynomial to be of the form $g = h(x^k)$ it is necessary that
the linear term (the trace) vanish. The set of traceless matrices is a
proper subvariety of $\Sp(2g),$ so by Theorem \ref{langweil}, for a
large $p,$ the probability that a given matrix is traceless is $\ll
1/p,$ and in particular is bounded away from $1.$ The proof of Theorem
\ref{randprodthm} now goes through verbatim to show that the set of
traceless matrices is asymptotically negligible in $\Sp(2g,
\mathbb{Z}).$

 Finally, the set of
cyclotomic polynomials of a given degree $2g$ is bounded by some
function $h(g).$ 
\end{proof}

\section{Free Group Automorphisms}
\label{free}

An automorphism of $\phi$ of  a free group $F_n$ is called
\emph{strongly irreducible}\footnote{This terminology, with strong
  support from this author, has been introduced by L. Mosher and
  M. Handel for what was previously known as \emph{irreducible with irreducible powers}} if no (positive) power of
$\phi$ sends a free factor $H$ of $F_n$ to a conjugate. This concept was
introduced by M.~Bestvina and M.~Handel \cite{bestvinahandel1}, and many of the results of the
theory of automorphisms of free groups are shown for such
automorphisms (for a survey, and the relationship between strongly
irreducible automorphisms and pseudo-Anosov automorphisms of surfaces,
the reader is urged to read the very clear survey\footnote{
to appear, but available at 
\newline
http://www.math.cornell.edu/~vogtmann/papers/AutQuestions/Questions.html}\cite{bridsonvogtmann}. By passing to the action of $\phi$ on the
abelianization of $F_n$ (equivalently, on $H_1(F_n, \mathbb{Z})$), Section
\ref{strong}\footnote{We need to change $\SL(n, \mathbb{Z})$ to
  $\GL(n, \mathbb{Z})$ throughout} shows the following:

\begin{theorem}
\label{irredpow}
Let $f_1, \dots, f_k$ be a generating set of the automorphism group of
$F_n.$ Consider all words of length $L$ in $f_1, \dots, f_k.$ Then,
for any $n,$ the probability that such a word is irreducible tends to
$1$ as $L$ tends to infinity and also the probability that
such a word is strongly irreducible tends to $1$ as $L$
tends to infinity.
\end{theorem}

\section{Galois groups of generic restricted polynomials}
\label{galoisrestricted}

Let $P_{N, d}(\mathbb{Z})$ be the set of monic polynomials of degree $d$ with
integral coefficients bounded by $N$ in absolute value. It is a
classical result of B.~L.~van~der~Waerden that the probability that
the Galois group of $p \in P_{N, d}(\mathbb{Z})$ is the full symmetric
group $S_d$ tends to $1$ as $N$ tends to infinity. The argument is
quite elegant: First, it is observed that a subgroup $H < S_d$ is the
full symmetric group if and only if $H$ intersects every conjugacy
class of $S_d.$ This means that $H$ has an element with every possible
cycle type. It is further noted that there is a cycle type $(n_1,
\dots, n_k)$ in the Galois group of $p$ over $\mathbb{Z}/p \mathbb{Z}$
if and only if $p$ factors over $\mathbb{Z}/p\mathbb{Z}$ into
irreducible polynomials of degrees $n_1, \dotsc, n_k.$ Using
Dedekind's generating function for the number of irreducible
polynomials over $\mathbb{Z}/p\mathbb{Z}$ of a given degree, it is
shown that the probability of a fixed partition is is bounded below by
a constant (independent of the prime $p$), and the proof is finished
by an application of a Chinese Remainder Theorem.

In this note, we ask the following simple-sounding question: Let $P_{N, d, a,
k}(\mathbb{Z})$   be the set of all polynomials in $P_{N,
d}(\mathbb{Z})$ where the coefficient of $x^k$ equals $a.$ Is it still
true that the Galois group of a random such polynomial is the full
symmetric group? The result would obviously follow if the probability
that the Galois group of a random general polynomial is ``generic''
were to go to $1$ sufficiently fast with $N.$ In fact, the probability
that an element of $P_{N, d}$ is \emph{reducible} (which means that
its Galois group is not transitive, hence not $S_n$) is of the order
of $1/N,$ so that approach does not work. 

Mimicking the proof of van der Waerden's result does not appear to work
(at least not easily): Dedekind's argument enumerates all irreducible
polynomials, and the result is not ``graded'' by specific
coefficients. It is certainly possible that the argument can be pushed
through, but this appears to be somewhat involved. 

Given this sad state of affairs, we first use a simple trick and Dirichlet's
theorem on primes on arithmetic progressions to show first the following technical result:
\begin{theorem}
\label{premainthm}
The probability that a random element of $P_{N, d, a, k}(\mathbb{Z}/p
\mathbb{Z})$ has a a prescribed splitting type $s$ 
approachs  the probability that a random unrestricted polynomial of degree
$d$ has the splitting type $s,$ as long as $p-1$ is relatively prime
to $(d-k)!,$ and as $p$ becomes large.
\end{theorem}
which implies (by van der Waerden's argument):

\begin{theorem}
\label{mainthm}
The probability that a random element of $P_{N, d, a, k}(\mathbb{Z})$
has $S_d$ as the Galois group tends to $1$ as $N$ tends to infinity,
\end{theorem}

It should be noted that the (multivariate) Large Sieve (as used by  P. X. Gallagher
in \cite{galgal}) can be used to give an effective estimate on the
probability in the statement of Theorems \ref{mainthm}: that
is: $p(N) \ll N^{-1/2} \log N.$

\subsection{Proof of Theorem \ref{mainthm}}
\label{mainproof}
We will need two ingredients other than van der Waerden's original 
idea. The first is A. Weil's estimate (Theorem \ref{weilcourbesthm}),
the second is the following classical result:
\begin{theorem}[\cite{lang}[Theorem VIII.9.1]]
\label{langthm}
Let $k$ be a field, and $n\geq 2$ an integer. Let $a\in k,$ $a\neq 0.$
Assume that for all prime numbers $p$ such that $p|n$ we have $a
\notin k^p,$ and if $4|n,$ then $a\notin - 4 k^4.$ Then $X^n - a$ is
irreducible in $k[X].$
\end{theorem}

Theorem \ref{langthm} goes essentially back to N. H. Abel's
foundational memoir.

We will need an additional observation:

\begin{lemma}
\label{freelem}
Let $q = p^l,$ and 
let $x_1, \dots, x_k \in \mathbb{F}_q.$ Let $a, b \in \mathbb{F}_p,$
with $(a, b) \neq (1, 0)$
and let $g(a, b)(x) = a x + b$ be a transformation of $\mathbb{F}_q$
to itself. Then, it is not possible for $g(a, b)$ to permute $x_1,
\dots, x_k,$ if $k!$ is coprime to $p-1.$
\end{lemma}

\begin{lemma}
\label{freecor}
Consider a polynomial $f$ of degree $d$ over $\mathbb{F}_p,$ such that
$d < p,$ and such that the coefficient of $x^{d-1}$ does not vanish. Then there is no pair $(a, b) \neq (1, 0),$ such that
$f(a x + b) = a^d f(x),$ for all $x \in \mathbb{F}_p.$
\end{lemma}
\begin{proof}
There are two distinct cases to analyze. The first is when $a=1.$ In
that case, $f(x+b) = f(x)$ for all $x\in \mathbb{F}_p,$ and since $p >
d,$ $f(x+b) = f(x),$ for all $x$ in the algebraic closure of $F_p.$ Let
$r$ be a root of $f.$ Then, so are $r+b, r+2b, \dotsc, r+b(p-1),$
which are all distinct as long as $b \neq 0,$ but
since $p$ is greater than $d$ that means that $f$ is identically $0.$

The second case is when $a\neq 1.$ In that case, $x_0 = b/(1-a)$ is
fixed under the substitution $x \rightarrow a x + b,$ so setting
$z = x - x_0,$ we see that $f(ax+b) = f(az),$ and so $f(ax+b) = a^d
f(x)$ implies that $f(az) = a^d f(z),$ for all $z \in \mathbb{F}_p.$
Since $p > \deg f,$ 
the corresponding
coefficients of the right and the left hand polynomials must be equal  Since the coefficient of
$x^{d-1}$ does not vanish, it follows tha $a=1,$ which contradicts our
assumption. 
\end{proof}

The argument now proceeds as follows. First, we note that if the
polynomial $f(x)$ of degree $d$ has a certain splitting type (hence Galois group)
over $\mathbb{F}_p$ then so does $f(a x + b)/a^d,$ for any $a\neq 0, b \in
\mathbb{F}_p.$ The set of all linear substitutions forms a group $\mathbb{A}$,
which acts freely on the set of polynomials of degree $d,$ except for
the (small) exceptional set of polynomials with a vanishing
coefficient of $x^{d-1}$ as long as $d < p$ (by Lemma \ref{freecor}), so the distribution of splitting types among the $\mathbb{A}$
orbits is the same as among all of the polynomials of degree $d.$ Now,
consider polynomials with constant term $1.$ How many of them are
there in the $\mathbb{A}$ orbit of $f(x)?$ It is easy to see that the
number is equal to the number of solutions to 
\[f(b) = a^d.\] If the curve $\mathcal{C}_f$ given by $f(x) - y^d$ is absolutely irreducible,
that number is $p + O(\sqrt{p}),$ by Theorem \ref{weilcourbesthm}. By
Theorem \ref{langthm}, in order for $\mathcal{C}_f$ to not be
absolutely irreducible, we must either have that $f(x) = g^q(x),$ for
some $q|d,$ or $f(x) = - 4 h^4(x),$ in case $4|d.$ But the number of
such polynomials is bounded by $O(p^{d/2}),$ which is asymptotically
neglible. So, we see that the distribution of splitting types amongst
polynomials of degree $d$ with constant term $1$ is the same as for
all polynomials, as long as $d < p.$ 

\section{Random walks on groups and graphs}
\label{intp2}

We consider the following situation:
$G$ is a finite \emph{undirected} (multi)graph with $n$ vertices, and
$\Gamma$ is a finite group. Each vertex $v_i$ of $G$ is decorated with
an element $t_i$ of $\Gamma;$ we assume that the set $\{t_1,
\dotsc, t_n\}$ generates $\Gamma.$ We make the following assumptions:
\begin{condition}
\label{standardcond}
\begin{enumerate}
\item The graph $G$ is
\emph{ergodic}: the adjacency matrix$A(G)$ 
of $G$ is a Perron-Frobenius matrix (meaning that there is a unique
eigenvalue of maximal modulus, and that some power $A^k(G)$ has all
entries positive). 
\item We assume further that for every nontrivial
one-dimensional unitary representation $\rho$ of $\Gamma,$ 
there exists $1\leq i, j\leq n,$ such that $\rho(t_i) \neq \rho(t_j).$
\end{enumerate}
\end{condition}
To each walk $w=(v_{i_1}, v_{i_2}, \dotsc, v_{i_k})$ on $G$ we
associate the element $\gamma_w = t_{i_k} \dots  t_{i_2} t_{i_1}.$ We
now denote the set of all walks of length $N$ by $W_N.$ 

Furthermore, we define a \emph{probability distribution} $P_N$ on
$\Gamma,$ as follows: the probability density $p_N: \Gamma \rightarrow
\reals_+$  is defined as $p_N(\gamma) = |\{w\in W_N \quad | \quad
\gamma_w = \gamma\}|/|W_N|.$

A slightly more abstract way to think of this is as follows: $P_N$ is
the function on the group ring $\reals[\Gamma],$ defined by:
\[
P_N = \dfrac{1}{|W_N} \sum_{w \in W_N} \gamma_w.
\]

Our main result is the following:
\begin{theorem}
\label{mythm1}
The distributions $P_N$ converge to the uniform distribution on
$\Gamma$ as $N$ goes to infinity. Furthermore, there is a constant $c
= c(G, \Gamma, T) < 1,$ such that $p_N(\gamma) - 1/|\Gamma| < c^N,$
for all $\gamma \in \Gamma.$ Here, $T$ refers to the assignment of
generators of $\Gamma$ to the vertices of $G.$
\end{theorem}

We will actually define a slightly stronger result: we will pick
$1\leq i, j \leq n,$ and consider the set of all walks $W_{N, i, j}$
of length $N$ from the $i$th to the $j$th vertex of $G.$ We define
distributions $P_{N, i, j}$ in the obvious way. Then:
\begin{theorem}
\label{mythm2}
Theorem \ref{mythm1} holds with $P_N$ replaced by $P_{N, i, j}.$
\end{theorem}

The starting point for the proof of the theorems above is Fourier
Transform on finite groups, which is discussed in Section
\ref{fouriergroups}. In particular, we will be using Theorem
\ref{fourierest} and Corollary \ref{closetoconst} to reduce the
question of whether a probability distribution is close to uniform to
the proving that the Fourier Transform is small at every
\emph{non-trivial} representation. The reader might well wonder how
moving the problem to Fourier transform space helps us -- the answer
is that it turns out that we can reduce the estimation of the
``fourier coefficients'' to questions in linear algebra, through the
construction in Section \ref{fouriertomat}.

\section{Fourier Transform on finite groups}
\label{fouriergroups}
For a thorough introduction to the topic of this section the reader is
referred to \cite{serrereps,simonreps}.
Let $\Gamma$ be a finite group, and let $f:\Gamma\rightarrow \cx$ be a
function on $\Gamma.$ Furthermore, let $\widehat{\Gamma}$ be the
\emph{unitary dual} of $\Gamma:$ the set of all irreducible complex unitary
representations of $\Gamma.$ To $f$ we can associate its \emph{Fourier
  Transform} $\hat{f}.$ This is a function which associates to
each $d$-dimensional unitary representation $\rho$ a $d\times d$ matrix
$\hat{f}(\rho)$ as follows:
\[
\hat{f}(\rho) = \sum_{\gamma \in \Gamma} f(\gamma) \rho(\gamma).
\]

There is an inverse transformation, as well. Given a function $g$
on $\widehat{\Gamma}$ which associates to each $d$-dimensional
representation $\rho$ a $d\times d$ matrix $g(\rho),$ we can write:
\[
g^\sharp (\gamma) = \dfrac{1}{\abs{\Gamma}} \sum_{\rho\in
\widehat{\Gamma}} d_\rho \tr (g(\rho) \rho(\gamma^{-1}),
\]
where $d_\rho$ is the dimension of $\rho.$ We mean ``inverse'' in the
most direct way possible:
\[
\hat{f}^\sharp = f.
\]

The following result is classical (see, eg, \cite{simonreps}):
\begin{theorem}
\[
\sum_{\rho in \widehat{\Gamma}} d_\rho^2 = \abs{\Gamma},
\]
\end{theorem}
and, together with the Fourier inversion formula, implies
\begin{theorem}
\label{fourierest}
Let $g$ be a function on $\widehat{\Gamma},$ such that for every
\emph{nontrivial} $\rho \in \widehat{\Gamma},$ 
\[
\opnorm{g(\rho)} < \epsilon,
\]
where $\opnorm{\bullet}$ denotes the operator norm (see Section
\ref{matnorm}). 
Then, for any $\gamma_1, \gamma_2 \in \Gamma,$
\[
\abs{g^\sharp(\gamma_1) - g^\sharp(\gamma_2)}  < 2\epsilon.
\]
\end{theorem}
\begin{proof}
First, note that for the trivial representation $\rho_0,$ the quantity
\[
d_{\rho_0}g(\rho_0)\rho_0(\gamma) = g(\rho_0),
\]
so does not depend on $\gamma.$
By the Fourier inversion formula, then,
\[
\begin{split}
\abs{g^\sharp(\gamma_1) - g^\sharp(\gamma_2)} & = \\
\left\lvert \dfrac{1}{\abs{\Gamma}}\sum_{\substack{\rho \in \widehat{\Gamma}\\
       \rho \neq \rho_0}} d_\rho \tr(g(\rho) (\rho(\gamma_1) -
   \rho(\gamma_2)))\right\rvert & \leq \\
\sum_{i=1}^2 \left\lvert \dfrac{1}{\abs{\Gamma}}\sum_{\substack{\rho \in \widehat{\Gamma}\\
       \rho \neq \rho_0}} d_\rho \tr\left(g(\rho) \rho(\gamma_i) \right)
   \right\rvert &
\underset{\text{by
     Eq. \eqref{traceineq}}}{\leq} \\
\dfrac{2}{\abs{\Gamma}} \sum_{\rho \in \widehat{\Gamma}}  d_\rho^2
\opnorm{g(\rho)} & < 2\epsilon.
\end{split}
\]
\end{proof}
\begin{corollary}
\label{closetoconst}
Under the assumption of Theorem \ref{fourierest}, and assuming in
addition that $g$ is real valued, if 
\[\sum_{\gamma \in \Gamma} g(\gamma) = 1,
\]
then \[g(\gamma) - 1/|\Gamma| < 2 \epsilon\quad\forall \gamma \in
\Gamma.\]
Furthermore, if $\Omega \in \Gamma,$ 
\begin{equation}
\label{closetoconstset}
\left| \sum_{\gamma \in \Omega} g(\gamma) -
  \dfrac{\Omega}{|\Gamma|}\right| < 2 \epsilon |\Omega.|
\end{equation}
\end{corollary}
\begin{proof}
Without loss of generality, suppose that $g(\gamma) > 1/|\Gamma|.$
Then there is a $\gamma_2,$ such that $g(\gamma_2) < 1/|\Gamma|.$
Thus,
\[
g(\gamma) - 1/|\Gamma| < g(\gamma)-g(\gamma_2) < 2\epsilon.
\]
The estimate \eqref{closetoconstset} follows immediately by summing
over $\Omega.$
\end{proof}

\section{Fourier estimates via linear algebra}
\label{fouriertomat}
In order to prove Theorem \ref{mythm2}, we would like to use Theorem
\ref{fourierest}, and to show the equidistribution result, we would
need to show that for every \emph{nontrivial} irreducible
representation $\rho,$ 
\begin{equation}
\label{coeffdecay}
\lim_{N\rightarrow \infty}\dfrac{1}{|W_{N, i, j}|}\tr{\sum_{w\in W_{N, i,
    j}}  \rho(\gamma_w)} = 0.
\end{equation}

To demonstrate Eq. \eqref{coeffdecay}, suppose that $\rho$ is
$k$-dimensional, so acts on a $k$-dimensional Hilbert space $H_\rho=H.$
Let $Z = L^2(G)$ -- the space of complex-valued functions from $V(G)$
to $\cx,$ let $e_1, \dotsc, e_n$ be the standard basis of $Z,$ and let
$P_i$ be the orthogonal projection on the $i$-th coordinate space. We
introduce the matrix 
\[
U_\rho = \sum_{i=1}^n P_i \otimes \rho(t_i) = 
\begin{pmatrix}
\rho(t_1) & 0 & \dots & 0\\
0 & \rho(t_2) & \dots & 0\\
\hdotsfor[2]{4}\\
0 & 0 & \dots & \rho(t_n)
\end{pmatrix},
\]
and also the matrix $A_\rho = A(G) \otimes I_H,$ where $I_H$ is the
identity operator on $H.$ Both $U_\rho$ and $A_\rho$ act on $Z \otimes H.$
The following is immediate:
\begin{lemma}
\label{matprodlem}
Consider the matrix $(U_\rho A_\rho)^l,$ and think of it as an
$n\times n$ matrix of $k\times k$ blocks. Then the $ij$-th block
equals the sum over all paths $w$ of length $l$ beginning at $v_i$ and
ending of $v_j$ of $\rho(\gamma_w).$
\end{lemma}

Now, let $T_{ji}$ be the operator on $Z$ which maps $e_k$ to
$\delta_{kj} e_i.$

\begin{lemma}
\label{projlemma}
\[
\tr{\left[\left((T_{ji}^t P_j)\otimes I_H\right) (U_\rho \otimes A_\rho)^N (P_i\otimes I_H)\right]} = \tr{\sum_{w\in  W_{N, i, j}} \rho(\gamma_w)}
\]
\end{lemma}

\begin{proof}
The argument of trace on the left hand side simply extracts the $ij$-th $k\times k$ block from $(U_{\rho} \otimes A_\rho)^N.$
\end{proof}

By submulticativity of operator norm, we see that 
\[
\opnorm{(T_{ji}^t P_j)\otimes I_H (U_\rho \otimes A_\rho)^N P_i\otimes
  I_H} \leq \opnorm{ (U_\rho \otimes A_\rho)^N}, 
\]
and so proving Theorem \ref{mythm2} reduces (thanks to Theorem
\ref{fourierest}) to showing
\begin{theorem}
\label{fundcollapse}
\[
\lim_{N\rightarrow \infty} \dfrac{\opnorm{(U_\rho \otimes
    A_\rho)^N}}{|W_{N, i, j}|
%\opnorm{A_\rho^N}
} = 0,
\]
for any non-trivial $\rho.$
\end{theorem}
\begin{notation}
We will denote the spectral radius of an operator $A$ by $\specrad{A}.$
\end{notation}
Since $|W_{N, I, j}| \asymp \mathcal{R}^N(A(G)),$ and by Gelfand's Theorem
(Theorem \ref{gelfand}), 
\[
\lim_{N\rightarrow \infty} \|B^N\|^{1/N} = \specrad{B},
\]
for any matrix $B$ and any matrix norm $\|\bullet\|,$ Theorem
\ref{fundcollapse} is equivalent to the statement that the spectral
radius of $U_\rho \otimes A_\rho$ is smaller than that of $A(G).$

Theorem \ref{fundcollapse} is proved in Section \ref{myproof}. 

\subsection{Proof of  Theorem  \ref{fundcollapse}}
\label{myproof}
\begin{lemma}
\label{twistlemma}
Let $A$ be a bounded hermitian operator $A:H\rightarrow H,$ and $U:
H\rightarrow H$ a unitary operator on the same Hilbert space $H.$
Then the spectral radius of $U A$ is smaller than the spectral radius
of $A,$ and the inequality is strict unless an eigenvector of $A$ with
maximal eigenvalue is also an eigenvector of $U.$
\end{lemma}
\begin{proof}
The spectral radius of $UA$ does not exceed the operator norm of $UA,$
which is equal to the spectral radius of $A.$ Suppose that the two are
equal, so that there is a $v,$ such that $\norm{UA v} = \specrad{A}{v},$
and $v$ is an eigenvector of $UA.$
Since $U$ is unitary, $v$ must be an eigenvector of $A,$ and since it
is also an eigenvector of $UA,$ it must also be an eigenvector of $U.$
\end{proof}

In the case of interest to us, $\rho$ is a $k$-dimensional irreducible
representation of $\Gamma,$ $U = \diag(\rho(t_1), \dots,
\rho(t_n),$ while $A = A(G) \otimes I_k.$ We assume that $A(G)$ is an
irreducible matrix, so that there is a unique eigenvalue of modulus
$\specrad{A(G)},$ that eigenvalue $\lambda_{\max}$ 
(the \emph{Perron-Frobenius eigenvalue})
is positive, and it has a strictly positive eigenvector $v_{\max}.$ We
know that the spectral radius of $A$ equals the spectral radius of
$A(G),$ and the eigenspace of $\lambda_{\max}$ is the set of vectors of
the form $v_{\max} \otimes w,$ where $w$ is an arbitrary vector in
$\mathbb{C}^k.$ If $v_{\max} = (x_1, \dotsc, x_n),$ we can write
$v_{\max} \otimes w = (x_1 w, \dotsc, x_n w),$ and so $U(v_{\max} \otimes
w) = (x_1 \rho(t_1) w, \dotsc, x_n \rho(t_n) w).$ Since all of the
$x_i$ are nonzero, in order for the inequality in Lemma
\ref{twistlemma} to be nonstrict, we must have some $w$ for which
$\rho(t_i) w = c w$ (where the constant $c$ does \emph{not} depend on
$i.$) 
Since the elements $t_i$ generate $\Gamma,$ the existence of such a
$w$ contradicts the irreducibility of $\rho,$ \emph{unless} $\rho$ is
one dimensional. This proves Theorem \ref{fundcollapse}

\section{Some remarks on matrix norms}
\label{matnorm}
In this note we use a number of matrix norms, and it is useful to
summarize what they are, and some basic relationships and inequalities
satisfied by them. For an extensive discussion the reader is referred
to the classic \cite{hornjohnson}. All matrices are assumed square,
and $n\times n.$

A basic tool in the inequalities below is the \emph{singular value
decomposition} of a matrix $A.$ 
\begin{definition}
The singular values of $A$ are the non-negative
square roots of the eigenvalues of
$A A^*,$ where $A^*$ is the conjugate transpose of $A.$ 
\end{definition}
Since $AA^*$
is a positive semi-definite Hermitian matrix for any $A,$ the singular
values $\sigma_1 \overset{\text{def}}{=} \sigma_{\max} \geq \sigma_2
\geq \dots$ are non-negative real numbers. For a Hermitian $A,$ the
singular values are simply the absolute values of the eigenvalues of $A.$

The first matrix norm is the \emph{Frobenius norm}, denoted by
$\norm{\bullet}.$
This is defined as 
\[
\norm{A} = \sqrt{\tr{A A^*}} = \sqrt{\sum_i \sigma_i^2}.
\]
This is also the sum of the square moduli of the elements of $A.$

The next matrix norm is the \emph{operator norm}, $\opnorm{\bullet},$
defined as 
\[
\opnorm{A} = \max_{\norm{v} = 1} \norm{Av} = \sigma_{\max}
\]

Both the norms $\norm{\bullet}$ and $\opnorm{\bullet}$ are
\emph{submultiplicative} (submultiplicativity is part of the
definition of matrix norm: saying that the norm $\matnorm{\bullet}$ is
submultiplicative means that $\matnorm{AB} \leq \matnorm{A}
\matnorm{B}$.)

From the singular value interpretation\footnote{A celebrated result of
  John von Neumann states that \emph{any} unitarily invariant matrix
  norm is a symmetric guage on the space of singular values -
  \cite{vnguage}.} of the two matrix norms and the Cauchy-Schwartz
inequality we see immediately that 
\begin{equation}
\label{twonorms}
\norm{A}/\sqrt{n}\leq \opnorm{A}\leq \norm{A}
\end{equation}

We will also need the following simple inequalities:
\begin{lemma}
\label{ulemma}
Let $U$ be a unitary matrix:
\begin{equation}
\label{traceineq}
\abs{\tr AU} \leq \norm{A}\sqrt{n} \leq n \opnorm{A}.
\end{equation}
\end{lemma}
\begin{proof}
Since $U$ is unitary, $\norm{U}= \norm{U^t} = \sqrt{n}.$
So, by the Cauchy-Schwartz inequality, 
$tr A U \leq \norm{A} \norm{U} = \sqrt{n} \norm{U}.$
The second inequality follows from the inequality \eqref{twonorms}.
\end{proof}

The final (and deepest result) we will have the opportunity to use is:
\begin{theorem}[Gelfand]
\label{gelfand}
For any 
operator $M,$ the spectral radius $\specrad{M}$ and any matrix norm
$\matnorm{\bullet},$
\[
\specrad{M} = \lim_{k\rightarrow \infty} \matnorm{M^k}^{1/k},
\]
\end{theorem}

\appendix

\section{Symplectic matrices and David Kirby's Theorem}
\label{kirbysec}
Recall that the \emph{symplectic quadratic form} in $2n$ dimensions is
given by the 
matrix
\[
J = 
\begin{pmatrix}
0 & \mathbf{I}_n\\
-\mathbf{I}_n & 0
\end{pmatrix},
\]
where $\mathbf{I}_n$ is the $n\times n$ identity matrix. A $2n \times
2n$ matrix $M$ is called \emph{symplectic} if it preserves the 
symplectic form $J,$ that is
\begin{equation}
\label{sympeq}
M^t J M = J.
\end{equation}
The invariance relation Eq. \eqref{sympeq} can be written slightly
more explicitly if we write $M$ in $n\times n$ block form, as 
\[
\label{genmat}
M = 
\begin{pmatrix}
A & B\\
C & D
\end{pmatrix}
\]
In that case, Eq. \eqref{sympeq} is equivalent to the following
system:
\begin{gather}
\label{symp1}
A^t C = C^t A (= (A^t C)^t),\\
\label{symp2}
B^t D = D^t B (= (B^t D)^t),\\
\label{symp3}
A^t D - C^t B = \mathbf{I}_n.
\end{gather}

This simple observation is all that is needed to prove the following
\begin{theorem}[David Kirby, \cite{Kirby_symplecticchar}]
\label{kirby}
Let $\mathcal{R}$ be a commutative ring with $1,$ and let $p(x) \in
\mathcal{R}[x]$ be a reciprocal polynomial.
Then, there is a symplectic matrix $M$ with coefficients in
$\mathcal{R},$  such that 
\[
\det (x \mathbf{I}_{2n} - M) = p(x).
\]
\end{theorem}
\begin{proof}[Proof sketch] 
In the notation of Eq. \eqref{genmat}, let $A = \mathbf{0},$ let $\det
B = 1,$ and let $C = -(B^t)^{-1}.$ Then, in order for \[
M
= \begin{pmatrix} 
A & B \\
C & D
\end{pmatrix} =
\begin{pmatrix}
\mathbf{0} & B \\
-(B^t)^{-1} & D
\end{pmatrix}
\]
to be symplectic, it is necessary and sufficient that $B^t D$ be
symmetric (since Eq. \eqref{symp1} and Eq. \eqref{symp3} hold
automatically.

Now, let 
\[
B_{ij} = \begin{cases}
0 & i + j < n\\
1 & i + j = n\\
b_{i+j+1-n} & i + j > n,
\end{cases}
\]
and let $D = E + F,$
where
\[
E_{ij} = \begin{cases}
0 & \mbox{$| i - j| > 1$ or $i=j > 1.$}\\
1 & | i - j| = 1\\
-1 & i=j=1,
\end{cases}
\]
and
\[
F_{ij} = \begin{cases}
0 & j < n\\
f_i & \mbox{otherwise}
\end{cases}
\]
for some $b_2, \dotsc, b_n, f_1, \dotsc, f_n \in \mathcal{R}.$
It is not hard to check that this is always so if $n = 1,$ and when
$n>1,$ it is necessary and sufficient that 
\begin{gather}
b_2 = f_1 + 1,\\
b_i = f_{i-1} + \sum_{j=2}^{i-1} b_j f_{i-j}.
\end{gather}
Now, if $p(x) = (1+x^{2n}) + a_1 * (x + x^{2n-1}) + \dotso + a_n(x^n +
x^{n+1}),$
a computation shows that a matrix $M$ as above with $f_i = a_{i-1} -
a_i$ does the trick. We won't go through all the details, but the main
idea in computing the determinant comes from the following:
\[
\begin{pmatrix}
x\mathbf{I} & -B\\
B^{-1} & x I - D
\end{pmatrix}
\begin{pmatrix}
B & 0\\
x \mathbf{I} & B^{-1}
\end{pmatrix} =
\begin{pmatrix}
0 & -\mathbf{I}\\
(x^2+1)\mathbf{I} - x D & x B^{-1} - D B^{-1},
\end{pmatrix}
\]
so by taking determinants:
\[
\det(x \mathbf{I} - M) = \det((x^2+1) \mathbf{I} - x D).
\]
\end{proof}

\bibliographystyle{plain}
\bibliography{rivin}
\end{document}